\documentclass[11pt,a4paper]{article}
\usepackage{amsfonts,amsmath,amssymb}
\usepackage{mathrsfs}
\usepackage[latin1]{inputenc}
\newtheorem{theo}{Theorem}
\newenvironment{Theo}{\begin{theo}\slshape}{\end{theo}}
\newtheorem{Lemm}{Lemma}[section]

\newtheorem{rema}{Remark}

\newcommand{\p}{\partial}
\newcommand{\Dom}{\mathrm{Dom}}

\def\qed{\hfill$\square$\par \bigskip}
\newenvironment{Demo}[1]{{\bf Proof #1.~}}{\qed}
\newcommand{\R}{\mathbb{R}}

\newcommand{\g}{\mathrm{\mathbf{g}}}
\newcommand{\f}{\mathrm{\mathbf{f}}}
\newcommand{\h}{\mathrm{\mathbf{h}}}

\newcommand{\HH}{\mathcal{H}}
\newcommand{\VV}{\mathcal{V}}
\newcommand{\norm}[1]{\left\Vert#1\right\Vert}
\newcommand{\abs}[1]{\left\vert#1\right\vert}
\newcommand{\set}[1]{\left\{#1\right\}}
\newcommand{\para}[1]{\left(#1\right)}
\newcommand{\cro}[1]{\left[#1\right]}

\newcommand{\To}{\longrightarrow}
\newcommand{\dive}{\textrm{div\,}}
\newcommand{\curl}{\textrm{curl\,}}
\newcommand{\D}{\mathrm{D}}
\newcommand{\B}{\mathrm{B}}
\newcommand{\U}{\mathrm{U}}
\newcommand{\V}{\mathrm{V}}
\newcommand{\W}{\mathrm{W}}
\newcommand{\ZZ}{\mathfrak{Z}}
\newcommand{\Uu}{\mathrm{U}_\sharp}
\newcommand{\Vv}{\mathrm{V}_\sharp}
\newcommand{\Ww}{\mathrm{W}_\sharp}
\usepackage{times}

\begin{document}
\title{\bf \Large Inverse boundary value problem for the dynamical heterogeneous Maxwell system}
\author{
\small {\bf Mourad Bellassoued\footnote{University of  Carthage, Faculty of Sciences of Bizerte, Department of Mathematics, 7021 Jarzouna Bizerte,
 Tunisia: mourad.bellassoued@fsb.rnu.tn}, \,\, Michel Cristofol\footnote{Aix-Marseille Universit\'e, CNRS, LATP, UMR 7353, 13453 Marseille, France: cristo@cmi.univ-mrs.fr}\,\, and Eric Soccorsi\footnote{Aix Marseille Université, CNRS, CPT, UMR 7332, 13288 Marseille, France: soccorsi@cpt.univ-mrs.fr}}
}

\date{}

\maketitle

\begin{abstract}
We consider the inverse problem of determining the isotropic inhomogeneous electromagnetic coefficients of the non-stationary Maxwell equations in a bounded domain of $\R^3$, from a finite number of boundary measurements. Our main result is a H\"older stability estimate for the inverse problem, where the measurements are exerted only in some boundary components.  For it, we prove a global Carleman estimate for the heterogeneous Maxwell's system with boundary conditions.
\end{abstract}

\section{Introduction}
In this paper, we discuss the uniqueness and stability in determining the isotropic electromagnetic coefficients
of the dynamical Maxwell equations, by boundary measurement of their solution. More precisely, given a continuous medium with dielectric permittivity $\lambda^{-1}$ and magnetic permeability $\mu^{-1}$, occupying an open, bounded and simply connected domain $\Omega\subset\R^3$ with ${\cal C}^\infty$ boundary
$\Gamma=\p\Omega$, and $T>0$, we consider the
following problem for the linear system of Maxwell's equations
\begin{equation}\label{1.1}
\begin{array}{lll}
\D'-\curl(\mu \B)=0, & \textrm{in }\,Q := \Omega \times (-T,T),\cr
\B'+\curl(\lambda \D)=0, & \textrm{in }\,Q,\cr
\dive\D=\dive\B=0, & \textrm{in }\,Q,\cr
\D\times\nu=0,\quad \B\cdot\nu=0, & \textrm{on }\,\Sigma
:= \Gamma\times(-T,T),
\end{array}
\end{equation}
where the prime stands for the time derivative. Here the electric induction $\D$ and the magnetic field $\B$ are three-dimensional vector-valued functions of the time $t$ and the space variable $x=(x_1,x_2,x_3)$, and $\nu=\nu(x)$ denotes the unit outward normal vector to $\Gamma$.
Moreover we attach the following initial condition to (\ref{1.1}):
\begin{equation}\label{1.3}
\B(x,0)=\B_0(x),\quad \D(x,0)=\D_0(x),\quad x\in\Omega.
\end{equation}
Assume that $\mu$ and $\lambda$ are scalar functions in $\mathcal{C}^2(\overline{\Omega})$ obeying
\begin{equation}\label{1.2}
\mu(x) \geq \mu_0, \quad   \lambda(x) \geq \lambda_0,\qquad x\in\overline{\Omega},
\end{equation}
for some $\lambda_0>0$ and $\mu_0>0$.
Next, in view of deriving existence and uniqueness results for (\ref{1.1}), introduce the following functional space
$$H(\curl, \Omega) := \{ u \in L^2(\Omega)^3,\ \curl u \in L^2(\Omega)^3 \}, $$
and denote by $\gamma_{\tau}$ the unique linear continuous application from $H(\curl, \Omega)$ into $H^{-1 \slash 2}(\Gamma)^3$, satisfying $\gamma_{\tau} u = u \wedge \nu$ when $u \in \mathcal{C}_0^{\infty}(\overline{\Omega})^3$
(see \cite{[DL3]}[Chap. IX A, Theorem 2]). Then, putting
$$H_0(\curl,\Omega) := \{ u \in H(\curl,\Omega),\ \gamma_{\tau} = 0 \}, $$
we see that the operator $i A$, where
$$ A \Phi := \left( \begin{array}{cc} 0  & \curl ( \mu .) \\ - \curl (\lambda . ) & 0 \end{array} \right),\ \Phi=( \D ,  \B ) \in \Dom(A):=H_0(\curl;\Omega) \times H(\curl;\Omega), $$
is selfadjoint in $\HH :=L^2(\Omega)^3 \times L^2(\Omega)^3$, endowed with the scalar product
$$
\langle \Phi,\tilde{\Phi} \rangle_{\HH} := \langle \lambda \D , \tilde{\D} \rangle_{{\rm L}^2(\Omega)^3} + \langle \mu \B , \tilde{\B} \rangle_{{\rm L}^2(\Omega)^3},\ \Phi=(\D,\B) \in \HH,\ \tilde{\Phi}=(\tilde{\D},\tilde{\B}) \in \HH.
$$
Further, in light of the last line of (\ref{1.1}), set
$$ H(\dive 0,\Omega) :=\{ u \in L^2(\Omega)^3,\,\dive u=0 \}\ {\rm and}\
H_0(\dive 0,\Omega) :=\{ u \in H(\dive 0,\Omega),\ \gamma_n u = 0 \}, $$
where $\gamma_n$ is the unique linear continuous mapping from $H(\dive ,\Omega) :=\{ u \in L^2(\Omega)^3,\,\dive u \in L^2(\Omega) \}$ onto $H^{-1 \slash 2}(\Gamma)$, such that $\gamma_n u = u \cdot \nu$ when
$u \in \mathcal{C}_0^{\infty}(\overline{\Omega})$ (see  \cite{[DL3]}[Chap. IX A, Theorem 1]).
Since $\HH_0 := {\rm H}(\dive 0;\Omega) \times {\rm H}_0(\dive 0;\Omega)$ is a closed subspace of $\HH$ and that $\HH_0^{\perp} \subset \ker A$, the restriction
$$ A_0 \Phi = A_{\HH_0} \Phi := A \Phi,\ \Phi \in \Dom(A_0) = \Dom (A) \cap \HH_0 := \VV, $$
is, by Stone's Theorem \cite{[DL5]}[Chap. XVII A,\S 4, Theorem 3], the infinitesimal generator of a unitary group of class $\mathcal{C}^0$ in $\HH_0$. Thus, by rewriting (\ref{1.1})-(\ref{1.3}) into the equivalent form
$$
\left\{ \begin{array}{l} \Phi'=A_0 \Phi \\ \Phi(0)=\Phi_0, \end{array} \right.\ {\rm with}\ \Phi=(  \D ,  \B )^T\ {\rm and}\ \Phi_0=( \D_0 , \B_0)^T,
$$
we get that:
\begin{Lemm}
\label{lm-existence}
Given $(\D_0,\B_0)\in \VV$ there exists a unique strong solution $(\D,\B)$ to (\ref{1.1}) starting from
$(\D_0,\B_0)$ within the following class
\begin{equation}
\label{1.4}
(\D,\B)\in \mathcal{C}^0(\R; \VV) \cap \mathcal{C}^1(\R; \HH ).
\end{equation}
\end{Lemm}
Moreover it holds true from \cite{[DL3]}[Chap. IX A, Remark 1] that
$$\VV= H_{\tau,0}(\curl,\dive 0;\Omega) \times {\rm H}_{n,0}(\curl,\dive 0;\Omega), $$
where
$$H_{*,0}(\curl,\dive 0;\Omega)=\{ u \in H^1(\Omega)^3,\ \dive u = 0\ {\rm and}\ \gamma_* u = 0 \},\ *=\tau, n.$$
For further reference we notice from Lemma \ref{lm-existence} that the solution
$(\D,\B)$ to (\ref{1.1})-(\ref{1.3}) actually satisfies:
\begin{equation}
\label{regularite}
(\D,\B) \in \cap_{p=0}^m \mathcal{C}^p([-T,T]; \Dom(A_0^{m-p})\ {\rm provided}\ (\D_0,\B_0) \in \Dom (A_0^m)\ {\rm and}\ \lambda,\ \mu \in \mathcal{C}^m(\overline{\Omega}),\ m \geq 1.
\end{equation}

The main purpose of this paper is to study the inverse problem of determining the dielectric permittivity $\lambda^{-1}=\lambda^{-1}(x)$ and the magnetic permeability $\mu^{-1}=\mu^{-1}(x)$ from a finite number of observations on the boundary $\Gamma$ of the solution $(\B, \D)$ to (\ref{1.1}) which corresponds to a realistic physical approach. This is an important problem not only in electromagnetics (see \cite{[Sta]})
but also in the identification of cracks/flaws in conductors
(see \cite{[HR]}) or the localization of lightning discharges (see \cite{[PHT]}). On the other hand, we obtain a reconstruction result which involves only a finite number of measurements which is not the case in most of the existing results.

\subsection{Inverse problem}
For suitable $\B^k_0$, $\D_0^k$, $k=1,2$, we aim to determine $\lambda(x)$, $\mu(x)$,
$x\in\Omega$, from the observation of
$$
\B^k_\tau(x,t),\quad\D_\nu^k(x,t),\quad (x,t)\in\Sigma,\quad k=1,2,
$$
where $\B_\tau=\B-(\B \cdot \nu)\nu$ (resp. $\D_\nu=(\D \cdot \nu)\nu$) denotes the tangential (resp. normal) component of $B$ (resp. $D$).

Notice that only a finite number of measurements are needed in the formulation of this inverse problem. For an overview of inverse problems for the Maxwell system, see the monograph \cite{[Romanov1]} by Romanov and Kabanikhin. For actual examples of inverse problems for the dynamical Maxwell system involving infinitely many boundary
observations (this is the case when the identification of the electromagnetic coefficients is made from the Dirichlet-to-Neumann map), we refer to Beleshev and Isakov \cite{[BI]}, Caro \cite{[C]}, Caro, Ola and Salo \cite{[CSO]}, Kurylev, Lassas and Somersalo \cite{[KLMS]}, Ola, Paivarinta and Somersalo \cite{[OPS]} and Salo, Kenig and Uhlmann
\cite{[SKU]}. It turns out that a small number of uniqueness and stability results for the inverse problem of determining the electromagnetic parameters of the Maxwell system with a finite number of measurements are available, such as \cite{[LY1], [LY2]}. In both cases, their proof is based on the methodology
of \cite{[BK]} or \cite{[IY3]}, which is by means of a Carleman estimate.
\medskip

For the formulation with a finite number of observations, Bukhgeim
and Klibanov \cite{[BK]} proposed a remarkable method based on a
Carleman estimate and established the uniqueness for similar inverse
problems for scalar partial differential equations.  See also Bellassoued \cite{[Be5]}, \cite{[Be6]}, Bellassoued
and Yamamoto \cite{[BeYa]}, \cite{[BeYa2]}, A. Benabdallah, M. Cristofol,
P. Gaitan and M. Yamamoto \cite{[BCGY]}, Bukhgeim \cite{[B]}, Bukhgeim,
Cheng, Isakov and Yamamoto \cite{[BCIY]}, Cristofol and Roques \cite{[CR]},  Cristofol and Soccorsi \cite{[CS]}, Imanuvilov and Yamamoto \cite{[IY2]}-\cite{[IY3]}, Isakov \cite{[I2]}, Kha\u\i darov \cite{[KH1]}, Klibanov \cite{[K1]}, \cite{[KL]}, Klibanov and Timonov \cite{[KT]}, Klibanov and Yamamoto \cite{[KY]}, Li and Yamamoto \cite{[LY1]}-\cite{[LY2]}, Yamamoto
\cite{[Y]}.

\medskip

A Carleman estimate is an inequality for a solution to a partial
differential equation with weighted $L^2$-norm and is a strong tool
also for proving the uniqueness in the Cauchy problem or the unique
continuation for a partial differential equation with non-analytic
coefficients. Moreover Carleman estimates have been applied
essentially for estimating the energy (e.g., Kazemi and Klibanov
\cite{[KK]}).
\medskip

As a pioneering work concerning a Carleman estimate, we refer to
Carleman's paper \cite{[Carleman]} which proved what is later called
a Carleman estimate and applied it for proving the uniqueness in the
Cauchy problem for a two-dimensional elliptic equation. Since
\cite{[Carleman]}, the theory of Carleman estimates has been studied
extensively. We refer to a general theory by H\"ormander \cite{[H]}
in the case where the symbol of a partial differential equation is
isotropic and functions under consideration have compact supports
(that is, they and their derivatives of suitable orders vanish on
the boundary of a domain). Later Carleman estimates for functions
with compact supports have been obtained for partial differential
operators with anisotropic symbols by Isakov \cite{[I2]}. Carleman
estimates for functions without compact supports, see Imanuvilov \cite{[Ima2]}, Tataru \cite{[T]}.  As for a direct derivation of
pointwise Carleman estimates for hyperbolic equations which are
applicable to functions without compact supports, see Klibanov and
Timonov \cite{[KT]}, Lavrent'ev, Romanov and Shishat$\cdot$ski\u\i \,  \cite{[La]}.
\medskip

The Carleman estimate for the non-stationary Maxwell's system was obtained for functions with compact supports, by Eller, Isakov,
Nakamura and Tataru \cite{[EINT]}. Li and Yamamoto \cite{[LY1]}-\cite{[LY2]}, prove a Carleman estimate for two-dimensional Maxwell's equations in
isomagnetic anisotropic media for functions with compact supports. Lemmas \ref{L2.1} and \ref{L2.2} are our Carleman estimate for the Maxwell's system whose solutions have not necessarily compact supports.
\medskip

By the methodology by \cite{[BK]} or \cite{[IY3]} with such Carleman estimates, several uniqueness and stability results are available for the inverse problem for the Maxwell's system (\ref{1.1}). That is, in \cite{[LY1]}-\cite{[LY2]} Li and Yamomoto established the uniqueness in determining three coefficients, using finite number of measurements.
\medskip

Li and Yamamoto \cite{[LY2]}, consider nonstationary Maxwell's equations in an anisotropic medium in the $(x_1,x_2,x_3)$-space, where equations of the divergences of electric and magnetic flux densities are also unknown. Then they  discuss an inverse problem of determining the $x$ 3-independent components of the electric current density from observations on the plane $x_3=0$ over a time interval and prove conditional stability in the inverse problem provided the permittivity and the permeability are independent of $x_3$.
\medskip

In \cite{[NP1]}, S.Nicaise and C.Pignotti, consider the Heterogenous Maxwell's system defined in an open bounded domain. Under
checkable conditions on the coefficients of the principal part they proved a Carleman
type estimates where some weighted $H^1$-norm of solution is dominated by the $L^2$ norm
of the boundary traces $\p_\nu U$ and $U_t$, modulo an interior lower-order term. Once
homogeneous boundary conditions are imposed the lower-order term can be absorbed
by the standard unique continuation theorem. Unfortunately, to our knowledge,
these results may not be applied directly to the linearized inverse problem associated
to the original problem.

Our argument is based on a new Carleman estimate. In comparison with \cite{[LY1]} and \cite{[NP1]},
our Carleman estimate is advantageous in the following two points:
\begin{itemize}
  \item We show a Carleman estimate which holds over the whole domain $Q$. We need not
assume that the functions under consideration have compact supports and so ours
is different from the Carleman estimates presented in \cite{[LY1]}, and we can establish
a H\"older estimate.

  \item We do not need a priori any unique continuation property and compactness/uniqueness
argument to absorb the lower-order interior term. In our approach, we establish
a Carleman estimates for $H^1$-solutions of the hyperbolic equation with variable
coefficients. This is essential to the proof of our main result, because here our problem
is involved with a source term and we cannot use the standard compactness/
uniqueness argument as in \cite{[NP1]}.
\end{itemize}


\subsection{Notations and statement of the main result}
\label{sec-NMR}
In this subsection we introduce some notations used throughout this text and state the main result of this article. Pick $x_0\in\R^3\backslash\overline{\Omega}$, set $c(x)=\mu(x) \lambda(x)$ for $x \in \overline{\Omega}$, $c_0=\mu_0 \lambda_0$ where $\mu_0$ and $\lambda_0$ are the same as in (\ref{1.2}), and assume that the following condition
\begin{equation}\label{1.6}
\frac{3}{2}\abs{\nabla\log c(x)}\abs{x-x_0}\leq 1-\frac{\rho}{c_0}, \quad x\in\overline{\Omega},
\end{equation}
holds true for some $\rho \in (0,c_0)$. This purely technical condition was imposed by the method we use to solve the inverse problem under study, which is by means of the Carleman estimate stated in Lemma \ref{L2.2} for any weight function $\psi_0$ satisfying the two Assumptions (A1) and (A2). More precisely, in the particular case where
\begin{equation}
\label{wf}
\psi_0(x) := |x-x_0|^2,\ x \in \Omega,
\end{equation}
then (\ref{1.6}) arises from the classical pseudo-convexity condition expressed by (\ref{2.5}). The somehow non-natural condition (\ref{1.6}) is thus closely related to the peculiar expression (\ref{wf}) in the sense that another choice of $\psi_0$ fulfilling (A1) and (A2) may eventually lead to a completely different condition on $c(x)$.

Next, for $M_0>0$ and two given functions $\mu^\sharp,\lambda^\sharp \in\mathcal{C}^2(\omega)$, where
$\omega=\Omega\cap \mathcal{O}$ for some neighbourhood $\mathcal{O}$ of $\Gamma$ in $\R^3$, we define the admissible set of unknown coefficients $\mu$ and $\lambda$ as
\begin{equation}\label{1.7}
\Lambda_\omega(M_0)=\set{(\mu,\lambda)\ \textrm{obeying}\,(\ref{1.2})\ {\rm and}\ (\ref{1.6})\,;
\norm{(\mu,\lambda)}_{\mathcal{C}^2(\overline{\Omega})}\leq M_0\ \textrm{and}
\,(\mu,\lambda)=(\mu^\sharp,\lambda^\sharp)\ \textrm{in}\ \omega}.
\end{equation}
Further, the identification of $(\lambda,\mu)$ imposing, as will appear in the sequel, that $(B,D)$ be observed twice, we consider two sets of initial data $(D_0^k,B_0^k)$, $k=1,2$,
\begin{equation}\label{1.8}
\D_0^k(x)=\para{d_1^k(x),d_2^k(x),d_3^k(x)}^\top,\quad\B_0^k(x)=\para{b_1^k(x),b_2^k(x),b_3^k(x)}^\top,
\end{equation}
and define the $12\times 6$ matrix
\begin{equation}\label{1.9}
\mathcal{K}(x)=\left(
              \begin{array}{cccccc}
                e_1\times\B_0^1 &  e_2\times\B_0^1 &  e_3\times\B_0^1 & 0 & 0 & 0 \\
                0 & 0 & 0 &  e_1\times\D_0^1 &  e_2\times\D_0^1 &  e_3\times\D_0^1 \\
                 e_1\times\B_0^2 &  e_2\times\B_0^2 &  e_3\times\B_0^2 & 0 & 0 & 0 \\
                0 & 0 & 0 &  e_1\times\D_0^2 &  e_2\times\D_0^2 &  e_3\times\D_0^2 \\
              \end{array}
            \right),\quad x\in\Omega.
\end{equation}
We then write $(\B_i^k(x,t),\D_i^k(x,t))$ the solution to (\ref{1.1}) with initial data $(\B_0^k,\D_0^k)$, $k=1,2$, where $(\mu_i,\lambda_i)$, $i=1,2$, is substituted for $(\mu,\lambda)$.

Finally, noting
$\mathscr{H}(\Sigma)=H^3(-T,T;L^2(\Gamma))\cap H^2(-T,T;H^1(\Gamma))$ the Hilbert space
equipped with the norm
$$
\norm{u}^2_{\mathscr{H}(\Sigma)}=\norm{u}^2_{H^3(-T,T;L^2(\Gamma))}+\norm{u}^2_{H^2(-T,T;H^1(\Gamma))},\quad u\in \mathscr{H}(\Sigma),
$$
we now may state the main result of this paper as follows :

\begin{Theo}\label{T.1} Let $T>c_0^{-1 \slash 2} \max_{x\in\overline{\Omega}}\abs{x-x_0}$ and pick $(\B_0^k,\D_0^k)\in (H^2(\Omega)^3 \times H^2(\Omega)^3) \cap \VV$, $k=1,2$, in such a way that there exists a $6 \times 6$ minor $m(x)$ of
the matrix $\mathcal{K}(x)$ defined in \eqref{1.9}, obeying:
\begin{equation}
\label{1.9b}
m(x) \neq 0,\ x \in \overline{\Omega\backslash\omega}.
\end{equation}
Further, choose $(\mu_i,\lambda_i)\in\Lambda_\omega(M_0)$, $i=1,2$, so that
\begin{equation}\label{1.13}
\norm{\para{\B_i^k,\D_i^k}}_{\mathcal{C}^3(-T,T;W^{2,\infty}(\Omega))}\leq M,\ k=1,2,
\end{equation}
for some $M>0$. Then there are two constants $C>0$ and $\kappa\in (0,1)$, depending on $\Omega$, $\omega$, $T$, $M$ and $M_0$, such that we have:
$$
\norm{\mu_1-\mu_2}_{H^2(\Omega)}+\norm{\lambda_1-\lambda_2}_{H^2(\Omega)}\leq C\para{\sum_{k=1}^2\para{\norm{\para{\B_1^k-\B_2^k}_\tau}_{\mathscr{H}(\Sigma)}+ \norm{\para{\D_1^k-\D_2^k}_\nu}_{\mathscr{H}(\Sigma)}}}^\kappa.
$$
\end{Theo}
Notice that the condition (\ref{1.9b}), which is independent of the choice of the unknown coefficients $\mu$ and $\lambda$,
actually relates on the initial functions in (\ref{1.3}). Moreover this condition is stable with respect to perturbations in $\mathcal{C}^2$. Namely, if $(\B_0^k,\D_0^k)$ obeys (\ref{1.9b}) then this is the case for $(\tilde{\B}_0^k,\tilde{\D}_0^k)$ as well, provided  $\max_{k=1,2}\norm{(\B_0^k,\D_0^k)-(\tilde{\B}_0^k,\tilde{\D}_0^k)}_{\mathcal{C}(\Omega)}$ is sufficiently small.
\medskip
Furthermore there are actual choices of $(\B_0^k,\D_0^k)$, $k=1,2$, satisfying (\ref{1.9b}). This can be seen by taking
$$
\B_0^1(x)=e_1,\quad\D_0^1(x)=e_3,\quad\B_0^2(x)=e_2,\quad\D_0^2(x)=e_2,\quad x\in \overline{\Omega \setminus \omega}
$$
and selecting the $6\times6$ minor formed by rows $2,3,4,9,10$ and $12$. \\
Theorem \ref{T.1} asserts H\"older stability in determining the principal part within
the class defined by (\ref{1.7}), under the assumption (\ref{1.13}). Notice from (\ref{regularite}) that such a condition is
automatically fulfilled for $\lambda, \mu \in \mathcal{C}^7(\overline{\Omega})$ by chosing the initial data $(\B_0^k,\D_0^k)$, $k=1,2$, in $\Dom(A_0^7)$ (which is a dense in $\HH$).\\
The proof of Theorem \ref{T.1} is based on a Carleman estimate stated in Lemma \ref{L2.2} under the
conditions (\ref{1.2}) and (\ref{1.6}). Notice that (\ref{1.6}), which is essential to our argument, is much stronger than the usual uniform ellipticity condition.

The remainder of the paper is organized as follows: a Carleman estimate for the Maxwell system (\ref{1.1}) is established in Section 2, while Section 3 contains the proof of Theorem \ref{T.1}.

\section{Carleman estimate for Maxwell's system}
\setcounter{equation}{0}
As already mentioned, this section is devoted to the derivation of a global Carleman estimate for the Maxwell system (\ref{1.1}).

\subsection{The settings}
Let us consider the following second order hyperbolic operator
\begin{equation}\label{2.1}
Pu=\p_t^2u(x,t)-\dive(c(x)\nabla u)+\mathscr{R}_1(x,t;\p)u,\quad x\in\Omega,\quad t\in\R,
\end{equation}
where $\mathscr{R}_1$ is a first order partial operator with $L^\infty(\Omega\times\R)$ coefficients, and $c\in\mathcal{C}^2(\overline{\Omega})$ obeys
\begin{equation}
\label{2.1b}
c(x)\geq c_0,\ x\in\overline{\Omega},
\end{equation}
for some positive constant $c_0$. Putting
\begin{equation}
\label{2.2}
a(x,\xi)=c(x)\abs{\xi}^2,\quad x\in\Omega,\quad \xi\in\R^3,
\end{equation}
and recalling the definition of the Poisson bracket of two given symbols $p$ and $q$,
$$
\set{p,q}(x,\xi)=\frac{\p p}{\p \xi}\cdot\frac{\p q}{\p x} -\frac{\p
p}{\p x}\cdot\frac{\p q}{\p \xi}=\sum_{i=1}^n\para{\frac{\p p}{\p
\xi_i}\frac{\p q}{\p x_i}-\frac{\p p}{\p x_i}\frac{\p q}{\p \xi_i}},
$$
we introduce two assumptions.\\
{\bf Assumption (A1)}.\ There exists $\psi_0 \in \mathcal{C}^2(\overline{\Omega};\R_+^*)$ satisfying
\begin{equation}\label{2.4}
\set{a,\set{a,\psi_0}}(x,\xi)>0,\quad x\in\overline{\Omega},\quad
\xi\in \R^3 \backslash\set{0},
\end{equation}
where $a$ is given by (\ref{2.2}).\\
Since $\overline{\Omega}$ is compact and $a(x,\xi)$ is a homogenous
function with respect to $\xi$, it is clear that (\ref{2.4}) yields the existence of some constant $\varrho>0$ such that
we have:
\begin{equation}
\label{2.5}
\frac{1}{4}\set{a,\set{a,\psi_0}}(x,\xi)\ge 2\varrho c(x) |\,\xi\,|^2,\quad
x\in\overline{\Omega},\quad \xi\in \R^3 \backslash\set{0}.
\end{equation}
{\bf Assumption (A2)}.\
The function $\psi_0(x)$ has no critical points on
$\overline{\Omega}$:
$$
\min_{x\in\overline{\Omega}}\abs{\nabla\psi_0(x)}^2 > 0.
$$
Further, $\varrho$ being  the same as in (\ref{2.5}), fix $\delta>0$ and $\beta \in (0,\varrho)$, in such a way that,
upon eventually enlarging $T$, we have:
\begin{equation}
\label{2.8}
\beta T^2>\max_{x\in\overline{\Omega}}\psi_0(x)+\delta.
\end{equation}
Hence, picking $\beta_0>0$ and setting
\begin{equation}
\label{2.7}
\psi(x,t)=\psi_0(x)-\beta t^2+\beta_0,\ x \in \overline{\Omega},\ t \in [-T,T],
\end{equation}
so that
$$
\min_{x\in\overline{\Omega}}\psi(x,0)\geq\beta_0,
$$
we check out from (\ref{2.8}) that
\begin{equation}
\label{2.11}
\psi(x,\pm T)\leq \beta_0-\delta,\ x\in\overline{\Omega}.
\end{equation}
Notice from (\ref{2.7}) and (\ref{2.11}) that
\begin{equation}\label{2.10}
\max_{x\in\overline{\Omega}}\psi(x,t)\leq \beta_0-\frac{\delta}{2},\ \abs{t} \in ( T-2\epsilon, T],
\end{equation}
for some constant $\epsilon\in (0,T/2)$.\\
In view of (\ref{2.7}) we may now recall the following global Carleman estimate for second order scalar hyperbolic equations, with weight function $\varphi : \Omega \times \R \To \R$ defined as
\begin{equation}
\label{2.12}
\varphi(x,t) = e^{\gamma \psi(x,t)},\ x \in \overline{\Omega},\ t \in [-T,T],
\end{equation}
for some fixed $\gamma>0$.

\begin{Theo}\label{T.2}
Assume (A1)-(A2). Then there exist two constants $C_0>0$ and $s_0>0$ such that for every $s\geq s_0$ the following Carleman estimate
\begin{eqnarray}
& & C_0\int_{Q}\!e^{2s\varphi}s\para{\abs{\nabla v}^2+\abs{\p_t v}^2+s^2\abs{v}^2} dx dt \nonumber \\
& \leq & \int_Q\! e^{2s\varphi}\abs{Pv(x,t)}^2 dx dt
+\int_{\Sigma}\!\! s e^{2s\varphi}\para{\para{\abs{\nabla v}^2+\abs{\p_t v}^2}+s^2\abs{v}^2} d\sigma dt,
\label{2.13}
\end{eqnarray}
holds true whenever $v\in H^1(Q)$ verifies $\p_t^jv(\pm T,\cdot)=0$ for $j=0,1$, and the right hand side of (\ref{2.13}) is finite. Here $P$ is defined by (\ref{2.1})-(\ref{2.1b}) and $d\sigma$ denotes the volume form of $\Gamma$.
\end{Theo}
For the proof see Bellassoued and Yamamoto \cite{[BeYa3]}, where this result is obtained from a direct computation based on integration by parts.
\medskip

It is worth mentioning that there are actual examples of functions $\psi_0$ fulfilling (A1)-(A2), provided the conductivity function $c$ defined in (\ref{2.1}) verifies (\ref{1.6})
for some $x_0 \in \R^3 \backslash\overline{\Omega}$ and $\varrho \in (0,c_0)$, where $c_0$ is the constant defined in (\ref{2.1b}).
Indeed, by putting $\psi_0(x)=\abs{x-x_0}^2$ and recalling (\ref{2.2}), we get through an elementary computation that
$$
\frac{1}{4}\set{a,\set{a,\psi_0}}(x,\xi)=2c^2(x)\para{1-\frac{\nabla\,c\cdot(x-x_0)}{2c}}|\,\xi\,|^2+
2c(\nabla c\cdot\xi)(\xi\cdot(x-x_0)),
$$
so (\ref{1.6}) immediately yields
$$
\frac{1}{4}\set{a,\set{a,\psi_0}}(x,\xi) \geq 2\varrho\,c |\,\xi\,|^2.
$$
This entails (A1) by (\ref{2.1b}). Moreover (A2) is evidently true as well since $\nabla\psi_0(x)\neq 0$ for every $x\in\overline{\Omega}$.

\subsection{Decoupling of the system of equations}
\label{sec-dec}
Consider now the following Maxwell system
\begin{equation}\label{2.17}
\begin{array}{lll}
\U'-\curl(\mu_1 \V)=\f, & \textrm{in }\,Q,\cr
\V'+\curl(\lambda_1 \U)=\g, & \textrm{in }\,Q,\cr
\dive\,\U=\dive\, \V=0, & \textrm{in }\,Q,\cr
\U\times\nu=0,\quad \V\cdot\nu=0, & \textrm{on }\,\Sigma,
\end{array}
\end{equation}
where the source terms $\f, \g \in H^1(Q;\R^3)$ satisfy the boundary condition
\begin{equation}\label{2.18}
\f(x,t)=\g(x,t)=0,\ (x,t)\in \omega\times(-T,T).
\end{equation}
For further reference we recall from (\ref{2.10}) that
\begin{equation}
\label{2.19}
\begin{array}{lll}
\max_{x\in\overline{\Omega}}\varphi(x,t)\leq d_0:=e^{\gamma(\beta_0-\delta/2)},\ \abs{t} \in [T-2\epsilon,T),\cr
\min_{x\in\overline{\Omega}}\varphi(x,0)\geq d_1:=e^{\gamma\beta_0},
\end{array}
\end{equation}
and then state the main result of \S \ref{sec-dec}:
\begin{Lemm}
\label{L2.1}
Assume (A1)-(A2) and let $\h=(\f,\g) \in H^1(Q;\R^3)^2$ obey (\ref{2.18}). Then we may find two constants $C_1>0$ and $s_1>0$, for which the Carleman estimate
\begin{eqnarray}
C_1\int_{Q} e^{2s\varphi} s\para{\abs{\nabla_{x,t} \W}^2+ s^2\abs{\W}^2} dx dt
& \leq & \int_Q e^{2s\varphi}\para{\abs{\nabla_{x,t} \h}^2+\abs{\h}^2} dx dt +s^3e^{2d_0s}\norm{\W}^2_{H^1(Q)} \nonumber \\
& & +\int_{\Sigma} s e^{2s\varphi}\para{\abs{\nabla\W}^2+\abs{\W'}^2+s^2\abs{\W}^2} d\sigma dt, \label{2.20}
\end{eqnarray}
is true for any $\W=(\U,\V)$ solution to the Maxwell system (\ref{2.17}), whenever $s\geq s_1$.
\end{Lemm}
\begin{Demo}{}
The first step of the proof involves bringing (\ref{2.17}) into two independent systems of decoupled equations. Namely, by differentiating the first line in (\ref{2.17}) with respect to $t$, and then substituting $\g-\curl(\lambda_1\U)$ for $\V'$ in the obtained equality, we obtain that
$$
\U''+\curl\para{\mu_1\curl(\lambda_1\U)}=\f'+\curl(\mu_1\g),\ {\rm in}\ Q.
$$
This entails
$\U''+\curl\para{\mu_1\lambda_1\curl\U}+\curl\para{\mu_1\nabla\lambda_1\times \U}=\f'+\curl(\mu_1\g)$, and hence
$$\U''+\mu_1\lambda_1\curl\para{\curl \U}+\nabla(\mu_1\lambda_1)\times\curl \U+\curl\para{\mu_1\nabla\lambda_1\times \U}=\f'+\curl(\mu_1\g),\ {\rm in}\ Q.
$$
From this, the well-known identity $\curl \curl \U = - \vec{\Delta}\ \U + \nabla \dive \U$ and the third line of (\ref{2.17}) then follows that
\begin{equation}
\label{2.25}
\U''-\mu_1\lambda_1\vec{\Delta} (\U)+\mathscr{R}_1 \U=\f'+\curl (\mu_1 \g),\ {\rm in}\ Q,
\end{equation}
where $\mathscr{R}_1=\mathscr{R}_1(x,\p)$ is some first order operator with bounded coefficients in $\Omega$.

Arguing in a similar way, we find that
\begin{equation}
\label{2.26}
\V''- \mu_1\lambda_1\vec{\Delta} (\V)+\mathscr{S}_1 \V=\g'-\curl (\lambda_1 \f),\ {\rm in}\ Q,
\end{equation}
for another first order operator $\mathscr{S}_1=\mathscr{S}_1(x,\p)$ with bounded coefficients in $\Omega$.\\
Therefore, putting (\ref{2.17}) and (\ref{2.25})-(\ref{2.26}) together, we end up getting that any solution $\W=(\U,\V)$ to the Maxwell system (\ref{2.17}) verifies
\begin{equation}
\label{2.27}
\left\{\begin{array}{lll}
\U''-\mu_1\lambda_1\vec{\Delta}(\U)+\mathscr{R}_1\U=G_1, & \textrm{in}\,\, Q\cr
\U\times\nu=0,\quad \curl(\lambda_1\U)\cdot\nu=0, & \textrm{on}\,\, \Sigma,
\end{array}
\right.
{\rm and}\
\left\{\begin{array}{lll}
\V''-\mu_1\lambda_1\vec{\Delta}(\V)+\mathscr{S}_1\V=G_2, & \textrm{in}\,\, Q\cr
\V\cdot\nu=0,\quad \curl(\mu_1\V)\times\nu=0, &  \textrm{on}\,\, \Sigma,
\end{array}
\right.
\end{equation}
where $G_1=\f'+\curl (\mu_1 \g)$ and $G_2=\g'-\curl(\lambda_1 \f)$.

Further, consider a cut-off function $\eta \in\mathcal{C}^\infty(\R;[0,1])$ fulfilling
\begin{equation}
\label{2.28}
\eta(t)=\left\{\begin{array}{ll}
1 & {\rm if}\ \abs{t}<T-2\epsilon, \cr
0 & {\rm if}\ \abs{t}\geq T-\epsilon,
\end{array}
\right.
\end{equation}
where $\epsilon$ is the same as in (\ref{2.10}), and set
$$
\Uu=\eta\U,\quad \Vv=\eta\V,\quad K_1=\eta G_1+2\eta'\U'+\eta''\U,\quad K_2=\eta G_2+2\eta'\V'+\eta''\V,
$$
in such a way that we have
$$
\left\{\begin{array}{lll}
\Uu''-\mu_1\lambda_1\vec{\Delta}(\Uu)+\mathscr{R}_1\Uu=K_1, & \textrm{in}\,\, Q\cr
\Uu\times\nu=0,\quad \curl{(\lambda_1\Uu)}\cdot\nu=0, &  \textrm{on}\,\, \Sigma
\end{array}
\right.
{\rm and}\
\left\{\begin{array}{lll}
\Vv''-\mu_1\lambda_1\vec{\Delta}(\Vv)+\mathscr{S}_1\Vv=K_2, & \textrm{in}\,\, Q\cr
\Vv\cdot\nu=0,\quad \curl{(\mu_1\Vv)}\times\nu=0, &  \textrm{on}\,\, \Sigma,
\end{array}
\right.
$$
directly from (\ref{2.27}).
Then, each of the two above systems being a principally scalar hyperbolic system, it follows from the two identities
\begin{equation}
\label{2.29}
\Uu(\cdot,\pm T)=\Uu'(\cdot,\pm T)=0,\quad \Vv(\cdot,\pm T)=\Vv'(\cdot,\pm T)=0,
\end{equation}
and Theorem \ref{T.2}, that $\Ww=(\Uu,\Vv)$ obeys the Carleman estimate
\begin{eqnarray}
& & C_0 \int_{Q} e^{2s\varphi} s\para{\abs{\nabla_{x,t} \Ww}^2+s^2\abs{\Ww}^2} dx dt \nonumber  \\
& \leq & \sum_{j=1,2} \int_Q\ e^{2s\varphi}\abs{K_j(x,t)}^2 dx dt
+\int_{\Sigma} s e^{2s\varphi}\para{\abs{\nabla\Ww}^2+\abs{\Ww'}^2+s^2\abs{\Ww}^2} d\sigma dt,
\label{2.32}
\end{eqnarray}
for all $s\geq s_0$. Here we have used Theorem \ref{T.2} for the diagonal system $\Uu''-\mu_1\lambda_1\vec{\Delta}(\Uu)$ and that we can absorb the non-decoupled first term $\mathscr{R}_1\Uu$.\\
Moreover, as $\eta'$ and $\eta''$ both vanish in $(-T+2\epsilon,T-2\epsilon)$ by (\ref{2.19}), there is a constant $C>0$ such that
$$
\sum_{j=1,2} \int_Q e^{2s\varphi}\abs{K_j(x,t)}^2dxdt\leq C\int_Q\! e^{2s\varphi}\para{\abs{\h}^2+\abs{\nabla_{x,t} \h}^2}dxdt+e^{2d_0s}\norm{\W}^2_{H^1(Q)},
$$
according to (\ref{2.29})-(\ref{2.32}) and since
$$
\int_{Q} e^{2s\varphi} s\para{\abs{\nabla_{x,t} \W}^2+ s^2\abs{\W}^2} dx dt\leq C \int_{Q} e^{2s\varphi} s\para{\abs{\nabla_{x,t} \Ww}^2+s^2\abs{\Ww}^2} dx dt +s^3e^{2d_0s}\norm{\W}^2_{H^1(Q)}
$$
 we obtain the result.
\end{Demo}
\subsection{Reduction of the boundary terms}
\label{sec-red}
The method used to derive a global Carleman estimate for the solution to (\ref{1.1}) is to replace
the local boundary problem (\ref{1.1}) in $\Omega \times(-T,T)$ by an equivalent
one stated on the half space $\R^3_+\times(-T,T)$. This is possible since the boundary $\Gamma$ can be represented as the zero level set of some $\mathcal{C}^\infty$ function in $\R^3$. Namely, $\Gamma$ being a $\mathcal{C}^\infty$ surface, there exist $\theta \in \mathcal{C}^\infty(\R^3)$ and some neighbourhood $\mathcal{V}$ of $\Gamma$ in $\R^3$
such that $\Gamma=\{ x \in \mathcal{V},\ \theta(x)=0 \}$. We choose $\mathcal{V}$ so small that $\mathcal{V} \subset \mathcal{O}$, where $\mathcal{O}$ is defined in \S \ref{sec-NMR}, write $y=(y_1,y_2,y_3)=(y',y_3)$ the system of normal geodesic coordinates where $y'=(y_1,y_2)$ are orthogonal coordinates in $\Gamma$ and $y_3=\theta(x)$ is the normal coordinate, and call $x=\Phi(y)$, where $\Phi'(y)>0$ for all $y\in\widehat{\mathcal{V}}:=\Phi^{-1}(\mathcal{V})$, the corresponding coordinates mapping. As
$$
\widehat{\Gamma}:=\Phi^{-1}(\Gamma)=\set{y\in\widehat{\mathcal{V}};\,\,y_3=0}\subset\R^2,
$$
we may assume that $\widehat{\mathcal{V}}=\Phi^{-1}(\mathcal{V})$ is a cylinder of the form $\widehat{\Gamma}\times(-r,r)$ with $r>0$.\\
Further, the Euclidean metric in $\R^3$ inducing the Riemannian metric with diagonal tensor $g$,
$$
g(y)={}^t\Phi'(y)\Phi'(y)=\textrm{Diag}(g_1,g_2,g_3),\quad y\in \widehat{\mathcal{V}},
$$
we use the notations of \cite{[EY]} and note $\cro{\frac{1}{\sqrt{g_1}}\frac{\p}{\p y_1},\frac{1}{\sqrt{g_2}}\frac{\p}{\p y_2},\frac{1}{\sqrt{g_3}}\frac{\p}{\p y_3}}$ the orthonormal basis associated by $g$ to the differential basis of vector fields $\cro{\frac{\p}{\p y_1},\frac{\p}{\p y_2},\frac{\p}{\p y_3}}$.
For any vector field $X(x)$ expressed with respect to the Euclidian basis $\cro{\frac{\p}{\p x_1},\frac{\p}{\p x_2},\frac{\p}{\p x_3}}$ as $X(x)=\sum_{i=1}^3\alpha^i(x)\frac{\p}{\p x_i}$,
we have an alternative representation $\widehat{X}(y)$ with respect to the new basis vectors $\cro{\frac{1}{\sqrt{g_1}}\frac{\p}{\p y_1},\frac{1}{\sqrt{g_2}}\frac{\p}{\p y_2},\frac{1}{\sqrt{g_3}}\frac{\p}{\p y_3}}$, given by
\begin{equation}
\label{2.37}
\widehat{X}(y)=\sum_{i=1}^3\widehat{\alpha}^i(y)\frac{1}{\sqrt{g_i}}\frac{\p}{\p y_i},\quad \widehat{\alpha}(y)={}^t\Psi(y)\alpha(\Phi(y)),\quad \Psi(y)=\Phi'(y)g^{-1/2}(y),\quad y\in  \widehat{\mathcal{V}}.
\end{equation}
The divergence (resp. curl) operator of any vector field $\widehat{X}(y)=\sum_{i=1}^3\widehat{\alpha}^i(y)\frac{1}{\sqrt{g_i}}\frac{\p}{\p y_i}$ with respect to the local coordinates $(y_1,y_2,y_3)$ is denoted by $\dive_g \widehat{X}$ (resp. $\curl_g \widehat{X}$), and can be brought into the form
\begin{equation}
\label{2.38}
\dive_g\widehat{X}=\frac{1}{\sqrt{\det g}}\sum_{j=1}^3\frac{1}{\sqrt{g_j}}\frac{\p}{\p y_j}\para{\sqrt{\det g}\,\widehat{\alpha}^j(y)}=
\sum_{j=1}^3\frac{1}{\sqrt{g_j}}\frac{\p\widehat{\alpha}^j}{\p y_j}+A(y)\cdot\widehat{X}(y),
\end{equation}
where $A(y)$ is some three dimensional vector (resp.
\begin{equation}
\label{2.39}
\curl_g\widehat{X}=\frac{1}{2}\sum_{i,j=1}^3\para{\frac{1}{\sqrt{g_j}}\frac{\p\widehat{\alpha}^i}{\p y_j}-\frac{1}{\sqrt{g_i}}\frac{\p\widehat{\alpha}^j}{\p y_i}}\frac{1}{\sqrt{g_j}}\frac{\p}{\p y_j}\times\frac{1}{\sqrt{g_i}}\frac{\p}{\p y_i}  +M(y)\widehat{X}(y),
\end{equation}
where $M(y)$ is some matrix function), whereas the outward normal vector field to $\widehat{\Gamma}$ at $y\in \widehat{\Gamma}$ is given by
\begin{equation}
\label{2.40}
\widehat{\nu}(y)=-\frac{1}{\sqrt{g_3}}\frac{\p}{\p y_3}.
\end{equation}
In light of (\ref{2.37})-(\ref{2.40}), we find out by performing the change of variable $x=\Phi(y)$ in
(\ref{2.17}), the space variable $x$ being restricted to be in $\Omega \cap \mathcal{V} \subset \omega \cap \mathcal{V}$, that
\begin{equation}\label{2.41}
\begin{array}{lll}
\widehat{\U}'-\curl_g(\widehat{\mu}_1 \widehat{\V})=0, & \textrm{in }\,\,\widehat{\mathcal{V}} \times (-T,T),\cr
\widehat{\V}'+\curl_g(\widehat{\lambda}_1 \widehat{\U})=0, & \textrm{in }\,\,\widehat{\mathcal{V}} \times (-T,T),\cr
\dive_g\,\widehat{\U}=\dive_g\, \widehat{\V}=0, & \textrm{in }\,\,\widehat{\mathcal{V}} \times (-T,T),\cr
\widehat{\U}\times\widehat{\nu}=0,\quad \widehat{\V}\cdot\widehat{\nu}=0, & \textrm{on }\,\,\widehat{\Gamma}\times(-T,T),
\end{array}
\end{equation}
where we have set
$$
\widehat{\mu}_1(y)={}^t\Psi(y)\mu_1(\Phi(y))\Psi(y),\quad \widehat{\lambda}_1(y)={}^t\Psi(y)\lambda_1(\Phi(y))\Psi(y),\quad y\in \widehat{\mathcal{V}}\subset\Phi^{-1}(\omega).
$$
Here we used the identity $\h(\Phi(y),t)=(\f(\Phi(y),t),\g(\Phi(y),t))=0$ for $t \in (-T,T)$ and $y\in \widehat{\mathcal{V}}$, arising from (\ref{2.18}).
Further, noting $\widehat{\U}=(\widehat{u}_1,\widehat{u}_2,\widehat{u}_3)$ and $\widehat{\V}=(\widehat{v}_1,\widehat{v}_2,\widehat{v}_3)$, the last equation in (\ref{2.41}) reads
$\widehat{u}_1=\widehat{u}_2=\widehat{v}_3=0$ on $\widehat{\Gamma}$, so we find that
\begin{equation}\label{2.45}
\widehat{\U}_{\widehat{\tau}}=0,\quad \widehat{\V}_{\widehat{\nu}}=0,\quad \textrm{on}\,\, \widehat{\Gamma},
\end{equation}
where $\widehat{\U}_{\widehat{\tau}}$ (resp. $\widehat{\V}_{\widehat{\nu}}$) denotes the tangential (resp. normal) component of $\widehat{\U}$ (resp. $\widehat{\V}$).
From this, (\ref{2.38}) and the third line in (\ref{2.41}) then follows that
$$
\frac{1}{\sqrt{g_3}}\frac{\p \widehat{u}_3}{\p y_3}=-\para{\frac{1}{\sqrt{g_1}}\frac{\p \widehat{u}_1}{\p y_1}+\frac{1}{\sqrt{g_2}}\frac{\p \widehat{u}_2}{\p y_2}}-A(y)\cdot\widehat{\U}=-A(y)\cdot\widehat{\U},
$$
whence
\begin{equation}\label{2.48}
\nabla \widehat{u}_3=\frac{1}{\sqrt{g_1}}\frac{\p \widehat{u}_3}{\p y_1}\widehat{e}_1+\frac{1}{\sqrt{g_2}}\frac{\p \widehat{u}_3}{\p y_2}\widehat{e}_2-\para{A(y)\cdot\widehat{\U}}\widehat{e}_3,\ {\rm on}\ \widehat{\Gamma}.
\end{equation}
Further, as
$\nabla \widehat{u}_1=\frac{1}{\sqrt{g_3}} \frac{\p \widehat{u}_1}{\p y_3}\widehat{e}_3$ and $\nabla \widehat{u}_2=\frac{1}{\sqrt{g_3}}\frac{\p \widehat{u}_2}{\p y_3}\widehat{e}_3$ on $\widehat{\Gamma}$, according to (\ref{2.45}), (\ref{2.48}) then yields
\begin{equation}\label{2.49}
\abs{\nabla\widehat{\U}}^2=\sum_{j=1}^3\abs{\nabla \widehat{u}_j}^2\leq C \para{\abs{\frac{\p \widehat{u}_1}{\p y_3}}^2+\abs{\frac{\p \widehat{u}_2}{\p y_3}}^2+\abs{\frac{\p \widehat{u}_3}{\p y_1}}^2+ \abs{\frac{\p \widehat{u}_3}{\p y_2}}^2+\abs{\widehat{\U}_{\widehat{\nu}}}^2},\ {\rm on}\ \widehat{\Gamma},
\end{equation}
where, for the sake of notational simplicity, we shall use the generic constant $C>0$ in the remaining of \S \ref{sec-red}. On the other hand, since
\begin{equation}
\label{2.50}
\curl_g\widehat{\U}\times\widehat{\nu}=\para{\frac{1}{\sqrt{g_3}}\frac{\p\widehat{u}_1}{\p y_3}-\frac{1}{\sqrt{g_1}}\frac{\p\widehat{u}_3}{\p y_1}}\frac{1}{\sqrt{g_1g_3}}\frac{\p}{\p y_1}- \para{ \frac{1}{\sqrt{g_2}}\frac{\p\widehat{u}_3}{\p y_2}-\frac{1}{\sqrt{g_3}}\frac{\p\widehat{u}_2}{\p y_3}}\frac{1}{\sqrt{g_2g_3}}\frac{\p}{\p y_2}+M(y)\widehat{\U}_{\widehat{\nu}},
\end{equation}
by (\ref{2.39})-(\ref{2.40}), we have
$$
\abs{\frac{\p \widehat{u}_1}{\p y_3}}^2+\abs{\frac{\p \widehat{u}_2}{\p y_3}}^2\leq C \para{\abs{\curl_g\widehat{\U}}^2+\abs{\widehat{\U}_{\widehat{\nu}}}^2+ \abs{\frac{\p \widehat{u}_3}{\p y_1}}^2+ \abs{\frac{\p \widehat{u}_3}{\p y_2}}^2},\ {\rm on}\ \widehat{\Gamma}.
$$
In view of (\ref{2.49}), this entails
\begin{eqnarray*}
\abs{\nabla\widehat{\U}}^2 & \leq  & C \para{\abs{\curl_g\widehat{\U}}^2+\abs{\frac{\p \widehat{u}_3}{\p y_1}}^2+ \abs{\frac{\p \widehat{u}_3}{\p y_2}}^2+\abs{\widehat{\U}_{\widehat{\nu}}}^2} \\ & \leq & C \para{\abs{\curl_g\widehat{\U}}^2+\abs{\nabla_{\widehat{\tau}} \widehat{\U}_{\widehat{\nu}}}^2+\abs{\widehat{\U}_{\widehat{\nu}}}^2},\ {\rm on}\ \widehat{\Gamma}.
\end{eqnarray*}
As a consequence we have
$$ \abs{\nabla\widehat{\U}}^2 \leq C \para{\abs{\curl_g ( \hat{\lambda}_1 \widehat{\U} ) }^2+\abs{\nabla_{\widehat{\tau}}\widehat{\U}_{\widehat{\nu}}}^2+\abs{\widehat{\U}_{\widehat{\nu}}}^2}, $$
whence
\begin{equation}
\label{2.52}
\abs{\nabla\widehat{\U}}^2 \leq C \para{\abs{\widehat{\V}'}^2+\abs{\nabla_{\widehat{\tau}}\widehat{\U}_{\widehat{\nu}}}^2+\abs{\widehat{\U}_{\widehat{\nu}}}^2},\ {\rm on}\ \widehat{\Gamma},
\end{equation}
by the second line of (\ref{2.41}).

Similarly, as $\dive_g \widehat{\V}=0$ in $\widehat{\mathcal{V}}$ from the third line of (\ref{2.41}), we get from (\ref{2.38}) that
$$
\frac{1}{\sqrt{g_3}}\frac{\p \widehat{v}_3}{\p y_3}=-\para{\frac{1}{\sqrt{g_1}}\frac{\p \widehat{v}_1}{\p y_1}+\frac{1}{\sqrt{g_2}}\frac{\p \widehat{v}_2}{\p y_2}}-A(y)\cdot\widehat{\V},\ {\rm on}\ \widehat{\Gamma}.
$$
This, combined with (\ref{2.45}), yields
$$
\abs{\nabla \widehat{v}_3}^2\leq C \para{\abs{\frac{\p \widehat{v}_1}{\p y_1}}^2+\abs{\frac{\p \widehat{v}_2}{\p y_2}}^2+\abs{\widehat{\V}_\tau}^2},
$$
and consequently
\begin{equation}\label{2.55}
\abs{\nabla\widehat{\V}}^2=\sum_{j=1}^3\abs{\nabla \widehat{v}_j}^2\leq C \para{\abs{\nabla_{\widehat{\tau}} \widehat{\V}_\tau}^2+\abs{\widehat{\V}_{\widehat{\tau}}}^2+\abs{\frac{\p \widehat{v}_1}{\p y_3}}^2+\abs{\frac{\p \widehat{v}_2}{\p y_3}}}^2,\ {\rm on}\ \widehat{\Gamma}.
\end{equation}
Furthermore, in light of (\ref{2.45}) and (\ref{2.50}) where $\widehat{\V}$ (resp. $\widehat{v}_j$, $j=1,2,3$) is substituted for $\widehat{\U}$ (resp. $\widehat{u}_j$, $j=1,2,3$), we see that $\abs{\frac{\p \widehat{v}_1}{\p y_3}}^2+\abs{\frac{\p \widehat{v}_2}{\p y_3}}^2$ is upper bounded, up to some multiplicative constant, by $\abs{\curl_g \widehat{\V}}^2 + \abs{\widehat{\V}_{\widehat{\tau}}}^2$, and hence
by $\abs{\curl_g(\widehat{\mu}_1\widehat{\V})}^2+\abs{\widehat{\V}_{\widehat{\tau}}}^2$, on $\widehat{\Gamma}$.
From this and the first line of (\ref{2.41}) then follows that
$$
\abs{\frac{\p \widehat{v}_1}{\p y_3}}^2+\abs{\frac{\p \widehat{v}_2}{\p y_3}}^2 \leq C \para{\abs{\widehat{\U}'}^2+\abs{\widehat{\V}_{\widehat{\tau}}}^2},\ {\rm on}\ \widehat{\Gamma},
$$
so, we end up getting with the aid of (\ref{2.55}):
\begin{equation}
\label{2.58}
\abs{\nabla\widehat{\V}}^2\leq C \para{\abs{\widehat{\U}'}^2+\abs{\nabla_{\widehat {\tau}} \widehat{\V}_{\widehat{\tau}}}^2+\abs{\widehat{\V}_{\widehat{\tau}}}^2}\ {\rm on}\ \widehat{\Gamma}.
\end{equation}
Finally, putting  (\ref{2.20}), (\ref{2.52}) and (\ref{2.58}) together, we may state the main result of \S \ref{sec-red}:

\begin{Lemm}
\label{L2.2}
Assume (A1)-(A2) and put $\h=(\f,\g)$. Then there are two constants $C_2>0$ and $s_2>0$  such that the following Carleman estimate
\begin{eqnarray*}
& & C_2s \int_{Q}\!e^{2s\varphi} \para{s^2\abs{\W}^2+ \abs{\nabla_{x,t} \W}^2} dx dt \\
& \leq & \int_Q\! e^{2s\varphi}\para{\abs{\h}^2+\abs{\nabla_{x,t} \h}^2} dx dt
+\mathscr{B}_{s,\varphi}(\W)+s^3e^{2d_0s}\norm{\W}^2_{H^1(Q)},
\end{eqnarray*}
where
\begin{equation}\label{2.60}
\mathscr{B}_{s,\varphi}(\W)=\int_{\Sigma}\!\! se^{2s\varphi}\para{\abs{\nabla_\tau\V_\tau}^2+\abs{\nabla_\tau\U_\nu}^2+\abs{\U_\nu'}^2+\abs{\V_\tau'}^2+s^2(\abs{\U_\nu}^2+\abs{\V_\tau}^2)} d\sigma dt,
\end{equation}
holds true for every solution $\W=(\U,\V)$ to (\ref{2.17}), provided $s \geq s_2$.
\end{Lemm}
\section{Inverse problem}
\setcounter{equation}{0}
This section contains the proof of Theorem \ref{T.1}, which is divided into five steps. Firstly, the unknown parameters $\lambda$ and $\mu$ are brought to the source term of the linearized system associated to (\ref{1.1}), governing the variation induced on the solution to (\ref{1.1}) by perturbating the permittivity by
$\lambda$ and the permeability by $\mu$.
The second  step follows the idea of Bukhgeim and Klibanov presented in \cite{[BK]}, which is to differentiate the linearized system with respect to $t$ in order to move the unknown coefficients in the initial condition.
The next step is to bound the energy of this system at time $t=0$ with the aid of the Carleman inequality of Theorem \ref{T.2}.
The fourth step involves relating $\lambda$ and $\mu$ to the above mentioned estimate through the Carleman inequality for stationary $({\rm div},{\rm curl})$-systems, stated in Lemma \ref{L3.2}. This is rather technical and lengthy so we proceed in a succession of the two Lemmas \ref{L3.3} and \ref{L3.5}. The last step, detailed in \S \ref{sec-CP}, is to derive the desired result from the estimates established in Lemmas \ref{L3.3}-\ref{L3.5}.

In the remaining of this text, $x_0$ is a fixed point in $\R^3 \setminus \overline{\Omega}$, we choose as in (\ref{wf})
$\psi_0(x):=|x-x_0|^2$ for every $x \in \overline{\Omega}$, and
$\varphi_0(x):=\varphi(x,0)$,
where $\varphi$ denotes the function defined by (\ref{2.7}) and (\ref{2.12}).
Moreover, for the sake of notational simplicity, we shall use the generic constant $C>0$ in the various estimates
of \S \ref{sec-PE}-\ref{sec-CP}.

\subsection{Linearized inverse problem}
Given $(\mu_i,\lambda_i) \in \Lambda_\omega(M_0)$, $i=1,2$, and $(\B_0^k,\D_0^k)\in H^2(\Omega)^3 \times H^2(\Omega)^3$, $k=1,2$, we
consider the solution $(\B_i^k,\D_i^k)$ to the system (\ref{1.1}) where $(\lambda_i,\mu_i)$ is substituted for $(\lambda,\mu)$, with initial condition (\ref{1.3}) where $(\B_0,\D_0)=(\B_0^k,\D_0^k)$. Hence, putting
$$\mu=\mu_1-\mu_2,\ \lambda=\lambda_1-\lambda_2,$$
and setting
\begin{equation}
\label{3.5}
\f_k=\curl(\mu\B^k_2),\quad \g_k=-\curl(\lambda\D^k_2),
\end{equation}
we find by a straightforward computation that $\U_k = \D^k_1-\D^k_2$ and $\V_k=\B^k_1-\B^k_2$ satisfy the system
\begin{equation}\label{3.3}
\begin{array}{lll}
\U'_k-\curl(\mu_1 \V_k)=\f_k, & \textrm{in }\,Q,\cr
\V'_k+\curl(\lambda_1 \U_k)=\g_k, & \textrm{in }\,Q,\cr
\dive\U_k=\dive\V_k=0, & \textrm{in }\,Q,\cr
\U_k\times\nu=0,\quad \V_k\cdot\nu=0, & \textrm{on }\,\Sigma,
\end{array}
\end{equation}
with the initial data
\begin{equation}
\label{3.4}
\U_k(x,0)=0,\quad \V_k(x,0)=0.
\end{equation}

Further, by using the following notations
\begin{equation}
\label{3.6}
X_{k,j}(x,t)=\p_t^j X_k(x,t) \ {\rm for}\ X=\U,\V,\f,\g\ {\rm and}\ j \in \mathbb{N}^*,
\end{equation}
it turns out by differentiating (\ref{3.3}) $j$-times with respect to $t$ that
\begin{equation}\label{3.10}
\begin{array}{lll}
\U'_{k,j}-\curl(\mu_1 \V_{k,j})=\f_{k,j}, & \textrm{in }\,Q,\cr
\V'_{k,j}+\curl(\lambda_1 \U_{k,j})=\g_{k,j}, & \textrm{in }\,Q,\cr
\dive\U_{k,j}=\dive\V_{k,j}=0, & \textrm{in }\,Q,\cr
\U_{k,j}\times\nu=0,\quad \V_{k,j}\cdot\nu=0, & \textrm{on }\,\Sigma,
\end{array}
\end{equation}
and, due to (\ref{3.5})-(\ref{3.3}), that $\U_{k,1}$ and $\V_{k,1}$ satisfy the initial condition:
\begin{equation}
\label{3.8}
\U_{k,1}(x,0)=\curl(\mu\B^k_0),\quad \V_{k,1}(x,0)=-\curl(\lambda \D^k_0).
\end{equation}
As will appear in \S \ref{sec-PE}-\ref{sec-CP}, the main benefit of dealing with (\ref{3.5})-(\ref{3.8}) in the analysis of the inverse problem of determining $\lambda$ and $\mu$, is the presence of these two unknown coefficients in the initial condition (\ref{3.8}).

\subsection{Preliminary estimates}
\label{sec-PE}
Let $j$ and $k$ be in $\{1, 2\}$. As $\W_{k,j}=(\U_{k,j},\V_{k,j})$ is solution to (\ref{3.10}), we notice from
Lemma \ref{L2.2} that
\begin{multline}\label{3.11}
C_2 s \int_{Q}\!e^{2s\varphi}\para{s^2\abs{\W_{k,j}}^2+\abs{\nabla_{x,t}\W_{k,j}}^2}\!dx dt  \cr
 \leq \int_Q\! e^{2s\varphi}\para{\abs{\h_{k,j}}^2+\abs{\nabla_{x,t}\h_{k,j}}^2}dx dt
+\mathscr{B}_{s,\varphi}(\W_{k,j})+s^3e^{2d_0s}\norm{\W_{k,j}}^2_{H^1(Q)} := \ZZ_{k,j}(s),
\end{multline}
for every $s \geq s_2$, where $\h_{k,j}=(\f_{k,j},\g_{k,j})$ and $\mathscr{B}_{s,\varphi}$ is given by (\ref{2.60}).
\medskip

Moreover by the assumption (\ref{1.9b}) we can derive that
\begin{equation}\label{3.110}
\sum_{k=1}^2\abs{\B_0^k(x)}^2\geq c_*,\quad\textrm{and}\quad \sum_{k=1}^2\abs{\D_0^k(x)}^2\geq c_*,\quad x\in \overline{\Omega\backslash\omega}
\end{equation}
for some positive constant $c_*$. Indeed, if $\B_0^1(x)\neq 0$ or $\B_0^2(x)\neq 0$ for all $x\in \overline{\Omega\backslash\omega}$ then $\sum_{k=1}^2\abs{\B_0^k(x)}^2>0$ in the compact set $\overline{\Omega\backslash\omega}$. Now if there exist $x_1\in \overline{\Omega\backslash\omega}$ such that $\B_0^1(x_1)=0$ then by (\ref{1.9b}) we have $\B_0^2(x)\neq 0$ for any $x\in \overline{\Omega\backslash\omega}$.
\medskip

Further we recall from \cite{[EY]} the following Carleman estimate for stationary $({\rm div},{\rm curl})$-systems:
\begin{Lemm}
\label{L3.2}
There exist two positive constants $s_3$ and $C_3$ depending only on $\psi_0$ and $\Omega$, such that we have
\begin{equation}
\label{3.12}
C_3s\int_\Omega e^{2s\varphi_0} \abs{\mathbf{u}}^2dx\leq \int_\Omega e^{2s\varphi_0} \para{\abs{\curl \mathbf{u}}^2+\abs{\dive \mathbf{u}}^2}dx,
\end{equation}
for every $s\geq s_3$ and $\mathbf{u}\in H_0^1(\Omega)$.
\end{Lemm}

Prior to the proof of Lemma \ref{L3.3}, we establish the following technical result, needed in the derivation of (\ref{3.13})-(\ref{3.13b}).
\begin{Lemm}
\label{L3.1}
There exists a constant $s_*>0$ depending only on $T$ such that we have
$$
\int_\Omega |z(x,0)|^2dx\leq 2 \para{ s \int_Q |z(x,t)|^2 dx dt + s^{-1} \int_Q |z'(x,t)|^2dxdt},
$$
for all $s \geq s_*$ and $z\in H^1(-T,T;L^2(\Omega))$.
\end{Lemm}
\begin{Demo}{}
Let $\eta \in \mathcal{C}^\infty([-T,T];[0,1])$ fulfills (\ref{2.28}) for some fixed $\epsilon \in (0,T \slash 2)$.
Then, the following identity
\begin{eqnarray*}
& & \int_\Omega |z(x,0)|^2dx =
 \int_{-T}^{0}\frac{d}{dt}\para{\int_\Omega\eta^2(t)|z(x,t)|^2dx}dt \\
&=& 2 \Re{\para{\int_{-T}^{0}\int_\Omega\eta^2(t)\overline{z(x,t)}z'(x,t)dxdt}} +2\int_{-T}^{0}\int_\Omega \eta(t)\eta'(t)|z(x,t)|^2dxdt,
\end{eqnarray*}
holds true for every $z\in H^1(-T,T;L^2(\Omega))$. Applying Young's inequality, this entails
$$
\int_\Omega |z(x,0)|^2dx \leq (s+2\| \eta' \|_{\infty}) \int_Q |z(x,t)|^2 dx dt +s^{-1} \int_Q |z'(x,t)|^2 dxdt,
$$
for each $s >0$, so the result follows by taking $s_*=2\| \eta' \|_{\infty}$.
\end{Demo}
Having said that, we are now in position to prove the:
\begin{Lemm}
\label{L3.3}
There exist two constants $C_4>0$ and $s_4>0$ such that the following estimates
\begin{equation}
\label{3.13}
s^2 \para{C_4 s \int_\Omega e^{2s\varphi_0} \para{\abs{\mu}^2+\abs{\lambda}^2}dx - \int_\Omega e^{2s\varphi_0} \para{\abs{\nabla\mu}^2 +\abs{\nabla\lambda}^2}dx} \leq \sum_{k=1}^2\ZZ_{k,1}(s),
\end{equation}
and
\begin{equation}
\label{3.13b}
C_4 s \int_\Omega e^{2s\varphi_0} \para{\abs{\nabla \mu}^2+\abs{\nabla \lambda}^2} dx
-\sum_{\abs{\alpha}= 2}\int_\Omega e^{2s\varphi_0} \para{\abs{\p^\alpha\mu}^2+\abs{\p^\alpha\lambda}^2}dx
\leq \sum_{j,k=1}^2 \ZZ_{k,j}(s),
\end{equation}
hold true for $k=1,2$, and $s \geq s_4$.
\end{Lemm}
\begin{Demo}{}
By applying Lemma  \ref{L3.1} for $z=e^{s\varphi} \U_{k,1}$, we get that
$$
C s^2\int_\Omega e^{2s\varphi_0} \abs{\U_{k,1}(x,0)}^2dx\leq s \int_Q e^{2s\varphi} \para{s^2\abs{\U_{k,1}(x,t)^2}+\abs{\U_{k,2}(x,t)}^2}dxdt,
$$
provided $s$ is large enough.
In light of (\ref{3.8})-(\ref{3.11}), this entails
\begin{equation}
\label{3.15}
C s^2\int_\Omega e^{2s\varphi_0} \abs{\curl(\mu\B^k_0)}^2dx \leq \ZZ_{k,1}(s).
\end{equation}
Further, taking into account that $\mu B_0^k \in H_0^1(\Omega)^3$ since $\mu$ vanishes in $\omega$ and $\dive \B_0^k=0$, we have
$$
Cs^3\int_\Omega e^{2s\varphi_0} \abs{\mu\B^k_0}^2dx - s^2\int_\Omega e^{2s\varphi_0} \abs{\nabla\mu}^2dx \leq s^2 \int_\Omega e^{2s\varphi_0} \abs{\curl(\mu\B^k_0)}^2 dx,
$$
by (\ref{3.12}), whence
\begin{equation}
\label{3.17}
Cs^3\int_\Omega e^{2s\varphi_0} \abs{\mu\B^k_0}^2dx - s^2\int_\Omega e^{2s\varphi_0} \abs{\nabla\mu}^2dx
\leq \ZZ_{k,1}(s),
\end{equation}
from (\ref{3.15}).
Similarly, by arguing as above with $z=e^{s\varphi} \V_{k,1}$ instead of $e^{s\varphi} \U_{k,1}$, we find some constant $C>0$ for which
$Cs^3\int_\Omega e^{2s\varphi_0} \abs{\lambda\D^k_0}^2dx - s^2\int_\Omega e^{2s\varphi_0} \abs{\nabla \lambda}^2 dx$ can be made smaller than
the right hand side of (\ref{3.17}) by taking $s$ sufficiently large. This, (\ref{3.110})
and (\ref{3.17}) entails (\ref{3.13}).

We turn now to showing (\ref{3.13b}). To do that we apply Lemma \ref{L3.1} with $z=e^{s\varphi} \p_i\U_{k,1}$, $i=1,2,3$, getting
$$
C \int_\Omega e^{2s\varphi_0} \abs{\p_i U_{k,1}(x,0)}^2 dx \leq s \int_Q e^{2s\varphi} \para{\abs{\nabla \U_{k,1}(x,t)}^2+s^{-2} \abs{\nabla \U_{k,2}(x,t)}^2}dxdt,
$$
for $s$ large enough. This yields
\begin{equation}
\label{3.21}
C \int_\Omega e^{2s\varphi_0} \abs{\p_i\curl(\mu\B^k_0)}^2dx \leq \sum_{j=1}^2 \ZZ_{k,j}(s),
\end{equation}
by (\ref{3.8})-(\ref{3.11}). Further, bearing in mind that $\dive \B_0^k=0$ and using that
$(\p_i \mu) B_0^k \in H_0^1(\Omega)^3$,
we obtain
\begin{equation}
\label{3.23}
C s \int_\Omega e^{2s\varphi_0} \abs{(\p_i\mu)\B^k_0}^2dx - \sum_{\abs{\alpha}= 2}\int_\Omega e^{2s\varphi_0} \abs{\p^\alpha \mu}^2dx \leq \int_\Omega e^{2s\varphi_0} \abs{\curl((\p_i \mu) \B^k_0)}^2dx,
\end{equation}
by (\ref{3.12}). Moreover, as $\curl ( ( \p_i \mu) B_0^k ) = \p_i \curl(\mu \B_0^k) - \mu \curl (\p_i B_0^k) - \nabla \mu \times \p_i B_0^k$ and $\mu \in H_0^1(\Omega)$, we have by applying the Poincare inequality
\begin{equation}
\label{3.22}
C\int_\Omega e^{2s\varphi_0} \abs{\curl((\p_i\mu) \B^k_0)}^2dx
\leq \int_\Omega e^{2s\varphi_0} \abs{\p_i\curl(\mu\B^k_0)}^2 dx + \int_\Omega e^{2s\varphi_0} \abs{\nabla \mu}^2 dx,
\end{equation}
hence
\begin{equation}
\label{3.22b}
C s \int_\Omega e^{2s\varphi_0} \abs{(\p_i\mu)\B^k_0}^2dx - \sum_{1 \leq |\alpha | \leq 2} \int_\Omega e^{2s\varphi_0} \abs{\p^{\alpha} \mu}^2 dx
\leq
\int_\Omega e^{2s\varphi_0} \abs{\p_i\curl(\mu\B^k_0)}^2 dx,
\end{equation}
by substituting the right hand side of (\ref{3.22}) for $\int_\Omega e^{2s\varphi_0} \abs{\curl((\p_i\mu) \B^k_0)}^2dx$ in (\ref{3.23}). Putting (\ref{3.21}) and (\ref{3.22b}) together, and summing up the obtained inequality over $i=1,2,3$,
we end up getting that
\begin{equation}
\label{3.24}
C s \int_\Omega e^{2s\varphi_0} \abs{\nabla \mu}^2dx - \sum_{1 \leq \abs{\alpha}\leq 2} \int_\Omega
e^{2s\varphi_0} \abs{\p^\alpha \mu}^2 dx \leq \sum_{j,k=1}^2 \ZZ_{k,j}(s).
\end{equation}
Here we used again (\ref{3.110}).
Finally, by arguing as before with $z=e^{s\varphi} \p_i \V_{k,1}$, $i=1,2,3$, instead of $\p_i \U_{k,1}$, we get that
(\ref{3.24}) remains true with $\mu$ replaced by $\lambda$. This completes the proof of (\ref{3.13b}).
\end{Demo}
Finally, we establish the:
\begin{Lemm}
\label{L3.5}
There are two constants $C_5>0$ and $s_5>0$ such that we have
$$C_5 \sum_{\abs{\alpha}= 2}\int_\Omega e^{2s\varphi_0} \para{\abs{\p^\alpha\mu}^2+\abs{\p^\alpha\lambda}^2}dx
- \int_\Omega e^{2s\varphi_0} \para{\abs{\nabla \mu}^2+\abs{\nabla \lambda}^2}dx
\leq \sum_{j,k=1}^2  \ZZ_{k,j}(s), $$
for all $s\geq s_5$.
\end{Lemm}
\begin{Demo}{}
In light of the two following basic identities
$\U_{k,1}(x,0)=\mu\curl\B^k_0+\nabla \mu \times \B_0^k$ and
$\V_{k,1}(x,0)=-\mu\curl\D^k_0-\nabla \lambda \times\D_0^k$, $k=1,2$, arising from (\ref{3.8}),
we have
$$
\mathcal{K}(x)\left(
                \begin{array}{c}
                  \nabla\mu \\
                  \nabla\lambda \\
                \end{array}
              \right)
=
\mathcal{A}(x)
 \left(
   \begin{array}{c}
     \mu \\
     \lambda\\
   \end{array}
 \right) - \left(
  \begin{array}{c}
    \U_{1,1}(x,0) \\
    -\V_{1,1}(x,0) \\
    \U_{2,1}(x,0) \\
    -\V_{2,1}(x,0) \\
  \end{array}
\right),\ {\rm with}\ \mathcal{A}(x)
=\left(
   \begin{array}{cc}
     \curl\B_0^1 & 0 \\
     0 & \curl\D_0^1 \\
     \curl\B_0^2 & 0 \\
     0 & \curl\D_0^2 \\
   \end{array}
 \right),
$$
hence
$$
\mathcal{K}(x)\left(
                \begin{array}{c}
                  \nabla\p_i\mu \\
                  \nabla\p_i\lambda \\
                \end{array}
              \right)
=
\p_i\mathcal{A} \left(
   \begin{array}{c}
     \mu \\
     \lambda\\
   \end{array}
 \right)
+\mathcal{A} \left(
   \begin{array}{c}
     \p_i\mu \\
     \p_i\lambda\\
   \end{array}
 \right)
-\p_i\mathcal{K}\left(
                \begin{array}{c}
                  \nabla\mu \\
                  \nabla\lambda \\
                \end{array}
              \right)
-\left(
  \begin{array}{c}
    \p_i\U_{1,1}(x,0) \\
    -\p_i\V_{1,1}(x,0) \\
    \p_i\U_{2,1}(x,0) \\
    -\p_i\V_{2,1}(x,0) \\
  \end{array}
\right),\  $$
for every $i=1,2,3$. From this and (\ref{1.9b}) then follows that
\begin{equation}\label{3.31}
\sum_{\abs{\alpha}=2}\para{\abs{\p^\alpha\mu}^2+\abs{\p^\alpha\lambda}^2}\leq C\para{\sum_{k=1}^2\para{\abs{\nabla\U_{k,1}(x,0)}^2+\abs{\nabla\V_{k,1}(x,0)}^2}+\sum_{\abs{\alpha}\leq 1}\para{\abs{\p^{\alpha} \mu}^2+\abs{\p^{\alpha} \lambda}^2}}.
\end{equation}
Further, by multiplying (\ref{3.31}) by $e^{s \varphi_0}$, integrating over $\Omega$, and upper bounding
$\int_{\Omega} e^{2 s \varphi_0} | \nabla \W_{k,1}(x,0) |^2 dx$, with the aid of Lemma \ref{L3.1}, we find out
that
\begin{eqnarray}
& & C \sum_{\abs{\alpha}= 2}\int_\Omega e^{2s\varphi_0} \para{\abs{\p^\alpha\mu}^2+\abs{\p^\alpha\lambda}^2}dx
- \int_\Omega e^{2s\varphi_0} \para{\abs{\nabla\mu}^2+\abs{\nabla\lambda}^2}dx \nonumber \\
& \leq &  s \sum_{k=1,2} \para{\int_{Q} e^{2 s \varphi} | \nabla \W_{k,1}(x,t)|^2  dx dt +
s^{-2} \int_{Q} e^{2 s \varphi} | \nabla \W_{k,2}(x,t)|^2 dx dt}. \label{3.30}
\end{eqnarray}
Here we took advantage of the fact that both $\mu$ and $\lambda$ belong to $H_0^1(\Omega)$ in order to get rid of the integral
$\int_\Omega e^{2s\varphi_0} \para{\abs{\mu}^2+\abs{\lambda}^2}dx$ by applying
the Poincar\'e inequality.
Evidently the result now follows from (\ref{3.30}) and Lemma \ref{L2.2}.
\end{Demo}

\subsection{Completion of the proof of the main result}
\label{sec-CP}
In light of (\ref{3.13b}) and Lemma \ref{L3.5} we may find $C>0$ such that
\begin{equation}
\label{3.32}
C \sum_{\abs{\alpha}\leq 2}\int_\Omega e^{2s\varphi_0} \para{\abs{\p^\alpha\mu}^2+\abs{\p^\alpha\lambda}^2}dx\leq \sum_{j,k=1}^2\ZZ_{k,j}(s),
\end{equation}
upon taking $s$ sufficiently large. Moreover, due to (\ref{3.5}), we have
\begin{eqnarray*}
& & \int_Q\! e^{2s\varphi}\para{\abs{\h_{k,j}}^2+\abs{\nabla \h_{k,j}}^2}dxdt \\
& \leq & \norm{(\B_2^k,\D_2^k)}_{\mathcal{C}^3(-T,T;W^{2,\infty}(\Omega))}^2 \left( \sum_{\abs{\alpha}\leq 2} \int_Q\! e^{2s\varphi}\para{\abs{\p^\alpha\mu}^2 + \abs{\p^\alpha\lambda}}^2 dx dt \right),
\end{eqnarray*}
for every $(j,k) \in \{1,2\}^2$, from where we get
\begin{eqnarray}
& & \sum_{\abs{\alpha}\leq 2} \left( C \int_\Omega e^{2s\varphi_0} \para{\abs{\p^\alpha\mu}^2+\abs{\p^\alpha\lambda}^2}dx
- \int_Q\! e^{2s\varphi}\para{\abs{\p^\alpha\mu}^2+\abs{\p^\alpha\lambda}^2}dx dt \right) \nonumber \\
& \leq & \sum_{j,k=1}^2 \left( \mathscr{B}_{s,\varphi}(\W_{k,j})
+s^3e^{2d_0s} \norm{\W_{k,j}}^2_{H^1(Q)} \right), \label{3.36}
\end{eqnarray}
by combining (\ref{3.11}) and (\ref{3.32}). Further, by recalling (\ref{2.7}) and (\ref{2.12}) we see for each $\alpha \in \{1,2,3\}^2$ with $|\alpha|=2$, that
\begin{equation}
\label{3.36b}
\int_Q\! e^{2s\varphi}\para{\abs{\p^\alpha\mu}^2+\abs{\p^\alpha\lambda}^2}dxdt
= \int_\Omega\! e^{2s\varphi_0}\para{\abs{\p^\alpha\mu}^2+\abs{\p^\alpha\lambda}^2} \mathcal{G}_s(x) dx,
\end{equation}
where
\begin{equation}
\label{3.36c}
\mathcal{G}_s(x) := \int_{-T}^Te^{-2s(\varphi_0(x)-\varphi(x,t))}dt \leq  \int_{-T}^Te^{-2s(1-\sigma(t))}dt :=
\mathfrak{g}(s),\ \sigma(t):=e^{-\gamma \beta t^2}.
\end{equation}
As $\lim_{s \rightarrow +\infty} \mathfrak{g}(s)=0$ by Lebesgue's Theorem,
we thus obtain from (\ref{3.36})-(\ref{3.36c}) that
\begin{equation}
\label{3.37}
C \sum_{\abs{\alpha}\leq 2}\int_\Omega e^{2s\varphi_0} \para{\abs{\p^\alpha\mu}^2+\abs{\p^\alpha\lambda}^2}dx
\leq \sum_{j,k=1}^2 \left( \mathscr{B}_{s,\varphi}(\W_{k,j})
+s^3e^{2d_0s} \norm{\W_{k,j}}^2_{H^1(Q)} \right),
\end{equation}
for $s$ sufficiently large.

Furthermore, bearing in mind that $x_0 \in \mathbb{\R}^3 \setminus \overline{\Omega}$, we notice from (\ref{2.7}) and (\ref{2.12}) that
\begin{equation}
\label{3.37a}
\varphi_0(x) \geq \min_{x \in \overline{\Omega}} e^{\gamma(|x-x_0|^2+\beta_0)} \geq d_1 > d_0,\ x \in \overline{\Omega},
\end{equation}
where $d_0$ is defined in (\ref{2.19}). From this and (\ref{3.37}) then follows that
\begin{equation}
\label{3.37b}
C\sum_{\abs{\alpha}\leq 2}\int_\Omega\para{\abs{\p^\alpha\mu}^2+\abs{\p^\alpha\lambda}^2}dx\leq
e^{Cs}\sum_{j,k=1}^2  \mathscr{B}(\W_{k,j})
+s^3e^{-2(d_1-d_0)s}M,
\end{equation}
where $M$ is the same as in (\ref{1.13}), and
\begin{equation}
\label{3.37c}
\mathscr{B}(\W)=\int_{\Sigma}\para{\abs{\nabla_\tau\V_\tau}^2+\abs{\nabla_\tau\U_\nu}^2+\abs{\U_\nu'}^2+\abs{\V_\tau'}^2+(\abs{\U_\nu}^2+\abs{\V_\tau}^2)} d\sigma dt,\ \W=(\U,\V).
\end{equation}
In view of (\ref{3.37a}), (\ref{3.37b})-(\ref{3.37c}) then yields
$$
\sum_{\abs{\alpha}\leq 2}\int_\Omega\para{\abs{\p^\alpha\mu}^2+\abs{\p^\alpha\lambda}^2}dx\leq
C\para{\sum_{j,k=1}^2 \mathscr{B}(\W_{k,j})}^\kappa,
$$
for some $\kappa \in (0,1)$, proving Theorem \ref{T.1}.
\medskip

\noindent {\bf Acknowledgements}

\noindent Part of this work was done while the first author was visiting the Universit\'e d'Aix-Marseille. He gratefully acknowledges the hospitality of the Universit\'e d'Aix-Marseille.




\end{document}